\documentclass[12pt]{article}      
\usepackage{amssymb}
\usepackage{amscd}

\input xy       
\xyoption{all}      

\newtheorem{theorem}{Theorem}[section]
\newtheorem{prop}[theorem]{Proposition}
\newtheorem{defn}[theorem]{Definition}
\newtheorem{lemma}[theorem]{Lemma}
\newtheorem{coro}[theorem]{Corollary}
\newtheorem{prop-def}{Proposition-Definition}[section]

\newcommand{\nc}{\newcommand}

\nc{\bin}[2]{ (_{\stackrel{\scs{#1}}{\scs{#2}}})}  
\nc{\binc}[2]{(\!\! \begin{array}{c} \scs{#1}\\
    \scs{#2} \end{array}\!\!)}  
\nc{\bincc}[2]{  ( {\scs{#1} \atop
    \vspace{-1cm}\scs{#2}} )}  
\nc{\bs}{\bar{S}}
\nc{\la}{\longrightarrow}
\nc{\rar}{\rightarrow}
\nc{\dar}{\downarrow}
\nc{\dap}[1]{\downarrow \rlap{$\scriptstyle{#1}$}}
\nc{\defeq}{\stackrel{\rm def}{=}}
\nc{\dis}[1]{\displaystyle{#1}}
\nc{\dotcup}{\ \displaystyle{\bigcup^\bullet}\ }
\nc{\hcm}{\ \hat{,}\ }
\nc{\hts}{\hat{\otimes}}
\nc{\hcirc}{\hat{\circ}}
\nc{\lleft}{[}
\nc{\lright}{]}
\nc{\curlyl}{\left \{ \begin{array}{c} {} \\ {} \end{array}
    \right .  \!\!\!\!\!\!\!} 
\nc{\curlyr}{ \!\!\!\!\!\!\!
    \left . \begin{array}{c} {} \\ {} \end{array}
    \right \} }
\nc{\longmid}{\left | \begin{array}{c} {} \\ {} \end{array}
    \right . \!\!\!\!\!\!\!}
\nc{\ora}[1]{\stackrel{#1}{\rar}}
\nc{\ola}[1]{\stackrel{#1}{\la}}
\nc{\scs}[1]{\scriptstyle{#1}}
\nc{\mrm}[1]{{\rm #1}}
\nc{\dirlim}{\displaystyle{\lim_{\longrightarrow}}\,}
\nc{\invlim}{\displaystyle{\lim_{\longleftarrow}}\,}
\nc{\dislim}[1]{\displaystyle{\lim_{#1}}}
\nc{\colim}{\mrm{colim}}
\nc{\mvp}{\vspace{0.3cm}}
\nc{\tk}{^{(k)}}
\nc{\tp}{^\prime}
\nc{\ttp}{^{\prime\prime}}
\nc{\svp}{\vspace{2cm}}
\nc{\vp}{\vspace{8cm}}
\nc{\proofend}{$\blacksquare$ \vspace{0.3cm}}
\nc{\modg}[1]{\!<\!\!{#1}\!\!>}
\nc{\intg}[1]{F_C(#1)}
\nc{\lmodg}{\!<\!\!}
\nc{\rmodg}{\!\!>\!}
\nc{\cpi}{\widehat{\Pi}}
\nc{\sha}{{\mbox{\cyr X}}}  
\nc{\shpr}{\diamond}    
\nc{\labs}{\mid\!}
\nc{\rabs}{\!\mid}

\nc{\ann}{\mrm{ann}}
\nc{\Aut}{\mrm{Aut}}
\nc{\can}{\mrm{can}}
\nc{\Cont}{\mrm{Cont}}
\nc{\rchar}{\mrm{char}}
\nc{\cok}{\mrm{coker}}
\nc{\dtf}{{R-{\rm tf}}}
\nc{\dtor}{{R-{\rm tor}}}

\nc{\Div}{{\mrm Div}}
\nc{\End}{\mrm{End}}
\nc{\Ext}{\mrm{Ext}}
\nc{\Fil}{\mrm{Fil}}
\nc{\Fr}{\mrm{Fr}}
\nc{\Frob}{\mrm{Frob}}
\nc{\Gal}{\mrm{Gal}}
\nc{\GL}{\mrm{GL}}
\nc{\Hom}{\mrm{Hom}}
\nc{\hsr}{\mrm{H}}
\nc{\hpol}{\mrm{HP}}
\nc{\id}{\mrm{id}}
\nc{\im}{\mrm{im}}
\nc{\incl}{\mrm{incl}}
\nc{\length}{\mrm{length}}
\nc{\mchar}{\rm char}
\nc{\mpart}{\mrm{part}}
\nc{\ql}{{\QQ_\ell}}
\nc{\qp}{{\QQ_p}}
\nc{\rank}{\mrm{rank}}
\nc{\rcot}{\mrm{cot}}
\nc{\rdef}{\mrm{def}}
\nc{\rdiv}{{\rm div}}
\nc{\rtf}{{\rm tf}}
\nc{\rtor}{{\rm tor}}
\nc{\res}{\mrm{res}}
\nc{\SL}{\mrm{SL}}
\nc{\Spec}{\mrm{Spec}}
\nc{\tor}{\mrm{tor}}
\nc{\Tr}{\mrm{Tr}}
\nc{\tr}{\mrm{tr}}

\nc{\ab}{\mathbf{Ab}}
\nc{\Alg}{\mathbf{Alg}}
\nc{\Bax}{\mathbf{Bax}}
\nc{\bfk}{{\bf k}}
\nc{\bfone}{{\bf 1}}
\nc{\detail}{\marginpar{\bf More detail}
    \noindent{\bf Need more detail!}
    \svp}
\nc{\Diff}{\mathbf{Diff}}   
\nc{\gap}{\marginpar{\bf Incomplete}\noindent{\bf Incomplete!!}
    \svp}
\nc{\FMod}{\mathbf{FMod}}
\nc{\Int}{\mathbf{Int}}
\nc{\Mon}{\mathbf{Mon}}
\nc{\proof}{\noindent{\bf Proof: }}
\nc{\remarks}{\noindent{\bf Remarks: }}
\nc{\Rep}{\mathbf{Rep}}
\nc{\Rings}{\mathbf{Rings}}
\nc{\Sets}{\mathbf{Sets}}
\nc{\bill}[1]{\marginpar{\bf To Bill}\noindent{\bf To Bill:}
    {\tt #1}\\ }
\nc{\li}[1]{\marginpar{\bf To Li}\noindent{\bf To Li:}
    {\tt #1}\\ }

\nc{\BA}{{\Bbb A}}
\nc{\CC}{{\Bbb C}}
\nc{\DD}{{\Bbb D}}
\nc{\EE}{{\Bbb E}}
\nc{\FF}{{\Bbb F}}
\nc{\GG}{{\Bbb G}}
\nc{\HH}{{\Bbb H}}
\nc{\LL}{{\Bbb L}}
\nc{\NN}{{\Bbb N}}
\nc{\QQ}{{\Bbb Q}}
\nc{\RR}{{\Bbb R}}
\nc{\TT}{{\Bbb T}}
\nc{\VV}{{\Bbb V}}
\nc{\ZZ}{{\Bbb Z}}

\nc{\cala}{{\cal A}}
\nc{\calc}{{\cal C}}
\nc{\cald}{\mathcal{D}}
\nc{\cale}{{\cal E}}
\nc{\calf}{{\cal F}}
\nc{\calg}{{\cal G}}
\nc{\calh}{{\cal H}}
\nc{\cali}{{\cal I}}
\nc{\call}{{\cal L}}
\nc{\calm}{{\cal M}}
\nc{\caln}{{\cal N}}
\nc{\calo}{{\cal O}}
\nc{\calp}{{\cal P}}
\nc{\calr}{{\cal R}}
\nc{\cals}{{\cal S}}
\nc{\calt}{{\cal T}}
\nc{\calw}{{\cal W}}
\nc{\calx}{{\cal X}}
\nc{\CA}{\mathcal{A}}

\nc{\fraka}{{\frak a}}
\nc{\frakB}{{\frak B}}
\nc{\frakm}{{\frak m}}
\nc{\frakp}{{\frak p}}
\nc{\frakS}{{\frak S}}
\nc{\frakA}{{\frak A}}
\nc{\frakx}{{\frak x}}

\font\cyr=wncyr10

\title{ 
On Free Baxter Algebras: Completions and the
Internal Construction
\thanks{The first author is supported in part by NSF grant
    \#DMS 97-96122. 
    MSC Numbers: Primary 16A06, 47B99. 
    Secondary 13A99,13B35,16W99.}}
\author{Li Guo and 
William Keigher
\\
Department of Mathematics and Computer Science\\
Rutgers University\\
Newark, NJ 07102 \\
(liguo@newark.rutgers.edu)\\
(keigher@newark.rutgers.edu)
}
\date{}

\begin{document}
\maketitle

\setcounter{section}{0}

\section{Introduction}
In a previous paper~\cite{G-K1}, we gave an explicit construction of
a free Baxter algebra. This construction is called the shuffle Baxter
algebra since it is described in terms of mixable shuffles. 
In this paper and its sequel~\cite{G-K2}, we will continue 
the study of free Baxter algebras.

There are two goals of this paper.
The first goal is to extend the construction of shuffle Baxter
algebras to completions of Baxter algebras. This process is motivated
by a construction of Cartier~\cite{Ca} and is analogous
to the process of completing a polynomial algebra to obtain
a power series algebra. However, as we will see later,
unlike the close similarity of properties of a polynomial algebra
and a power series algebra, properties of
a shuffle Baxter algebra and its completion 
can be quite different.

The second goal is to establish a connection between the
shuffle Baxter algebra we have constructed
to the standard Baxter algebra
constructed by Rota~\cite{Rot1}.
The shuffle Baxter algebra is an external construction in the
sense that it is a free Baxter algebra obtained without reference
to any other Baxter algebra. On the other hand, the standard
Baxter algebra is an internal construction,
obtained as a Baxter subalgebra inside a naturally
defined Baxter algebra constructed originally by Baxter~\cite{Ba}.
There are several restrictions on Rota's original construction
of standard Baxter algebras.
By modifying Rota's
method and making use of the shuffle Baxter algebra construction,
we are able to construct the standard Baxter algebra in full
generality.
The shuffle product construction of a free Baxter algebra
has the advantage that its module
structure and Baxter operator can be easily described.
The description of a free Baxter algebra as a standard Baxter algebra
has the advantage that its multiplication is very simple. 
We will give an explicit description of the isomorphism
between the shuffle Baxter algebra and
the standard Baxter algebra. 
This description will enable us to make use of properties of both 
the shuffle product description and
Rota's description of free Baxter algebras.
Some applications will be given in~\cite{G-K2}.

We will start with a brief summary of definitions and basic
properties of the shuffle Baxter algebra in section~\ref{sec:back}.
We also take the opportunity to extend the construction of shuffle
Baxter algebras to the category of Baxter algebras not necessarily
having an identity. 
In section~\ref{sec:comp}, we define the completion of a Baxter algebra
by making use of the filtration given by the Baxter operator,
and give a description of a free complete Baxter algebra in terms
of mixable shuffles. 
In section~\ref{sec:rota}, we construct the standard 
Baxter algebras, generalizing Rota.
Variations of the construction for complete Baxter algebras
and for Baxter algebras not necessarily having an identity
are also considered.  

\section{Shuffle Baxter algebras}
\label{sec:back}

We write $\NN$ for the additive monoid of
natural numbers $\{0,1,2,\ldots\}$ and 
$\NN_+=\{ n\in \NN\mid n>0\}$ for the positive integers. 

Let $\Rings$ denote the category of commutative rings with identity.
For any $C\in \Rings$, let $\Alg_C$ denote the category of
$C$-algebras with identity. 
For $C\in \Rings$
and for any $C$-modules $M$ and $N$,
the tensor product $M\otimes N$ is taken over $C$
unless otherwise indicated. 
Let $M$ be a $C$-module. 
For $n\in \NN$, denote
\[ M^{\otimes n}=\underbrace{M\otimes\ldots\otimes M}
    _{n\ {\rm factors}}\]
with the convention that $M^{\otimes 0}=C$.

\subsection{Baxter algebras}
Baxter algebras were first studied by Baxter~\cite{Ba}
and the category of Baxter algebras was first studied by
Rota~\cite{Rot1}. 
We recall basic definitions and properties of Baxter algebras.
See~\cite{G-K1,Rot1} for details.

\begin{defn}
Let $C$ be a ring, $\lambda\in C$,
and let $R$ be a $C$-algebra.
\begin{itemize}
\item
A {\bf Baxter operator of weight $\lambda$ on $R$ over $C$}
is a $C$-module
endomorphism $P$ of $R$ satisfying 
\begin{equation}
 P(x)P(y)=P(xP(y))+P(yP(x))+\lambda P(xy),\ x,\ y\in R.
\label{eq:bax1}
\end{equation}
\item
    A {\bf Baxter C-algebra of weight $\lambda$}
is a pair $(R,P)$ where $R$ is a $C$-algebra
and $P$ is a Baxter operator of weight $\lambda$
on $R$ over $C$.
\item
Let $(R,P)$ and $(S,Q)$ be two Baxter $C$-algebras
of weight $\lambda$.
A {\bf homomorphism of Baxter $C$-algebras}
$f:(R,P)\rar (S,Q)$ is a homomorphism
$f: R \rar S$ of $C$-algebras with the property that
$ f(P(x))=Q(f(x))$
for all $x \in R$.
\end{itemize}
\end{defn}
If the meaning of $\lambda$ is clear, we will
suppress $\lambda$ from the notation. 
Note that our $\lambda$ is $-q$ in the notation of Rota~\cite{Rot3}.

Let $\Bax_{C,\lambda}$ denote the category of Baxter
$C$-algebras of weight $\lambda$. 
A {\bf Baxter ideal} of $(R,P)$ is an ideal $I$ of $R$ such
that $P(I)\subseteq I$.
Other concepts of $C$-algebras, such as subalgebra and
quotient algebra, can also be defined for Baxter algebras~\cite{G-K1}.

\subsection{Shuffle Baxter algebras with an identity}
\label{sec:shuf}
Let $A\in \Alg_C$. 
In a previous work~\cite{G-K1}, we used mixable shuffles to construct
a mixable shuffle algebra $\sha_C(A)$ and proved that
it is a free Baxter $C$-algebra on $A$.

For $m,n\in \NN_+$, 
define the set of {\bf $(m,n)$-shuffles} by 
\begin{eqnarray*}
\lefteqn{ S(m,n)=} \\
&& \left \{ \sigma\in S_{m+n}
    \begin{array}{ll} {} \\ {} \end{array} \right .
\left | 
\begin{array}{l}
\sigma^{-1}(1)<\sigma^{-1}(2)<\ldots<\sigma^{-1}(m),\\
\sigma^{-1}(m+1)<\sigma^{-1}(m+2)<\ldots<\sigma^{-1}(m+n)
\end{array}
\right \}.
\end{eqnarray*}
Given an $(m,n)$-shuffle $\sigma\in S(m,n)$, 
a pair of indices $(k, k+1)$,\ $1\leq k< m+n$ is 
called an {\bf admissible pair} for $\sigma$
if $\sigma(k)\leq m<\sigma(k+1)$. 
Denote $\calt^\sigma$ for the set of admissible pairs for $\sigma$. 
For a subset $T$ of $\calt^\sigma$, call the pair
$(\sigma,T)$ a {\bf mixable $(m,n)$-shuffle},
where $(\sigma,T)$ is identified with $\sigma$
if $T$ is the empty set. 
Denote
\[ \bs (m,n)=\{ (\sigma,T)\mid \sigma\in S(m,n),\
    T\subset \calt^\sigma\} \]
for the set of {\bf $(m,n)$-mixable shuffles}. 

For $m,n \in \NN_+$, 
denote
$x=x_1\otimes\ldots\otimes x_m\in A^{\otimes m}$,
and
$y=y_1\otimes \ldots\otimes y_n\in A^{\otimes n}$. 
For $\sigma\in S_m$, denote
\[\sigma (x)=
    x_{\sigma(1)}\otimes x_{\sigma(2)} \otimes
    \ldots \otimes x_{\sigma(m)}.\]
Denote
$x\otimes y =x_1\otimes \ldots \otimes x_m
    \otimes y_1 \otimes \ldots\otimes y_n
    \in A^{\otimes (m+n)}$. 
and, for $\sigma\in S_{m+n}$,
denote
\[\sigma (x\otimes y) =u_{\sigma(1)}\otimes u_{\sigma(2)} \otimes
    \ldots \otimes u_{\sigma(m+n)},\]
where
\[ u_k=\left \{ \begin{array}{ll}
    x_k,& 1\leq k\leq m,\\
    y_{k-m}, & m+1\leq k\leq m+n. \end{array}
    \right . \]

\begin{defn}
Let $x\in A^{\otimes m}$, $y\in A^{\otimes n}$
and $\sigma\in S(m,n)$.
\begin{enumerate}
\item
$\sigma (x\otimes y)\in A^{\otimes (m+n)}$ is called 
a {\bf shuffle} of $x$ and $y$.
\item
Let $T$ be a subset of $\calt_\sigma$. 
The element 
\[ \sigma(x\otimes y; T)= u_{\sigma(1)}\hts u_{\sigma(2)} \hts
    \ldots \hts u_{\sigma(m+n)}, \]
where for each pair $(k,k+1)$, $1\leq k< m+n$, 
\[ u_{\sigma(k)}\hts u_{\sigma(k+1)} =\left \{\begin{array}{ll}
    u_{\sigma(k)} u_{\sigma(k+1)},  &
     (k,k+1)\in T\\
    u_{\sigma(k)}\otimes u_{\sigma(k+1)}, &
    (k,k+1) \not\in T
    \end{array} \right . \]
is called a {\bf mixable shuffle} of $x$ and $y$.
\end{enumerate}
\end{defn}

Fix a $\lambda\in C$. 
For $x=x_0\otimes x_1\otimes\ldots \otimes x_m\in
A^{\otimes (m+1)}$ and
$y=y_0\otimes y_1\otimes\ldots\otimes y_n\in A^{\otimes (n+1)}$,  
define

\[x \shpr y\
=\sum_{(\sigma,T)\in \bs (m,n)} \lambda^{\mid T\mid }
    x_0y_0\otimes \sigma(x\otimes y;T) 
    \in \bigoplus_{k\leq m+n+1} A^{\otimes k}.\]
Then $\shpr $ extends to a mapping
\[ \shpr : A^{\otimes (m+1)}\times A^{\otimes (n+1)} \rar
    \bigoplus_{k\leq m+n+1} A^{\otimes k}, m, n\in \NN\]
by $C$-linearity. 
Let

\[ \sha_C(A)=\sha_C(A,\lambda)= \bigoplus_{k\in\NN}
    A^{\otimes (k+1)}
= A\oplus A^{\otimes 2}\oplus \ldots. \]
Extending by additivity, the map 
$\shpr$ gives a $C$-bilinear map 

\[ \shpr: \sha_C(A) \times \sha_C(A) \rar \sha_C(A) \]
with the convention that
\[ A\times A^{\otimes (m+1)} \rar A^{\otimes (m+1)} \]
is the scalar multiplication on the left $A$-module
$A^{\otimes (m+1)}$. 
Define a $C$-linear endomorphism $P_A$ on
$\sha_C(A)$ by assigning
\[ P_A( x_0\otimes x_1\otimes \ldots \otimes x_n)
=\bfone_A\otimes x_0\otimes x_1\otimes \ldots\otimes x_n, \]
for all
$x_0\otimes x_1\otimes \ldots\otimes x_n\in A^{\otimes (n+1)}$
and extending by additivity.
Let $j_A:A\rar \sha_C(A)$ be the canonical inclusion map. 
We proved the following theorem in~\cite{G-K1}. 

\begin{theorem}
\label{thm:shua}
\begin{enumerate}
\item
The $C$-module $\sha_C(A)$, together with the multiplication
$\shpr$, is a commutative $C$-algebra with an identity.
\item
$(\sha_C(A),P_A)$, together with the natural embedding
$j_A:A\rightarrow \sha_C(A)$,
is a free Baxter $C$-algebra on $A$ (of weight $\lambda$). 
In other words,
for any Baxter $C$-algebra $(R,P)$ and any
$C$-algebra map
$\varphi:A\rar R$, there exists
a unique Baxter $C$-algebra homomorphism
$\tilde{\varphi}:(\sha_C(A),P_A)\rar (R,P)$ such that
the  diagram
\[\xymatrix{
A \ar[rr]^(0.4){j_A} \ar[drr]_{\varphi}
    && \sha_C(A) \ar[d]^{\tilde{\varphi}} \\
&& R } \]
commutes. 
\end{enumerate}
\end{theorem}

The Baxter $C$-algebra $(\sha_C(A),P_A)$ will be
called the {\bf shuffle Baxter $C$-algebra
(of weight $\lambda$) on $A$.}
When there is no danger of confusion, we will often suppress
the symbol $\shpr$ and simply denote $x y$ for $x\shpr y$
in $\sha_C(A)$. 

For a given set $X$, 
let $C[X]$ be the polynomial $C$-algebra on $X$
with the natural embedding $X\hookrightarrow C[X]$. 
Let $(\sha_C(X),P_X)$ be the Baxter $C$-algebra
$(\sha_C(C[X]),P_{C[X]})$.
$(\sha_C(X),P_X)$ will be
called the {\bf shuffle Baxter $C$-algebra
(of weight $\lambda$) on $X$.}

\begin{prop}
\label{co:shux}
$(\sha_C(X),P_X)$, together with the set embedding
\[j_X: X\hookrightarrow C[X] \ola{j_{C[X]}} \sha_C(C[X]),\]
is a free Baxter $C$-algebra on the set $X$,
described by the following universal property:
For any Baxter $C$-algebra $(R,P)$ over $C$ and any set map
$\varphi: X\rar R$, there exists
a unique Baxter $C$-algebra homomorphism
$\tilde{\varphi}:(\sha_C(X),P_X)\rar (R,P)$ such that
the diagram 
\[\xymatrix{
X \ar[rr]^(0.4){j_X} \ar[drr]_{\varphi}
    && \sha_C(X) \ar[d]^{\tilde{\varphi}} \\
&& R } \]
commutes. 
\end{prop}

If we choose $A=C$ in the construction of free Baxter
$C$-algebras, then we get
\[ \sha_C(C)=\bigoplus_{n=0}^\infty C^{\otimes (n+1)}
= \bigoplus_{n=0}^\infty C \bfone^{\otimes (n+1)}, \]
where
$\bfone^{\otimes (n+1)}
= \underbrace{\bfone_C \otimes \ldots \otimes \bfone_C}
_{(n+1)-{\rm factors}}$.
Thus $\sha_C(C)$ is a free $C$-module on the basis
$\bfone^{\otimes n}, n\geq 1$.

\begin{prop}
\label{prop:unit}
For any $m,n\in \NN$,
\[ \bfone^{\otimes (m+1)} \shpr \bfone^{\otimes (n+1)} =
\sum_{k=0}^m \binc{m+n-k}{n}\binc{n}{k} \lambda^k
\bfone^{\otimes (m+n+1-k)}.\]
\end{prop}

\subsection{Shuffle Baxter algebras without an identity}
We now construct a shuffle Baxter algebra in the category
of Baxter algebras not necessarily having an identity.

For $C\in \Rings$, let
$\Alg_C^0$ be the category of $C$-algebras not necessarily having
an identity and let 
$\Bax_C^0$ be the category of Baxter $C$-algebras not necessarily
having an identity. 
For $A\in \Alg_C^0$, we will use mixable shuffles to construct
a free Baxter algebra on $A$ in $\Bax_C^0$. This construction
was given in a special case in~\cite{G-K1}.
A similar construction can be carried out if $C$ is assumed
to be a commutative ring not necessarily having an identity,
but we will not give details here. 

Let $C\in \Rings$ and $A\in \Alg_C^0$ be given.
We use a well-known construction~\cite{Ca,Ja1} to embed $A$ in an
element $A^+ \in \Alg_C$. 
Let $A^+= C \oplus A$ with the addition defined componentwise
and the multiplication defined by
\[ (c,a) (d,b) = (cd, cb+da+ab), c,\ d\in C,\ a,\ b\in A.\]
Then $A^+$ is in $\Alg_C$ with $(\bfone_C,0)$ as the identity and
$a\mapsto (0,a)$ embeds $A$ in $A^+$ as a subobject in $\Alg_C^0$
(in fact, as an ideal).

Define 
\[ \sha_C(A)^0 = \oplus_{n\in \NN} ((A^+)^{\otimes n}\otimes A)\]
with the convention that $(A^+)^{\otimes 0}=C$.
Thus $\sha_C(A)^0$ is the $C$-submodule of $\sha_C(A^+)$
generated by tensors of the form
\[ x_0 \otimes\ldots \otimes x_n,\
x_i\in A^+, 0\leq i\leq n-1,\ x_n\in A.\]
Since any mixable shuffle of $x_0\otimes\ldots\otimes x_m$
and $y_0\otimes\ldots\otimes y_n$ has 
either $x_m$ or $y_n$ or $x_my_n$ as the last tensor factor, 
we see that $\sha_C(A)^0$ is a $C$-subalgebra of
$\sha_C(A^+)$. It is also clearly closed under the Baxter
operator $P_{A^+}$.
So $\sha_C(A)^0$, with the restriction of $P_{A^+}$,
denoted by $P_A$, is a subobject of $\sha_C(A^+)$ in $\Bax_C^0$.
It is called {\bf the shuffle Baxter algebra} on $A$ (of weight
$\lambda$) in the category $\Bax_C^0$. 

\begin{prop}
\label{prop:nonu}
$(\sha_C(A)^0,P_A)$, together with the natural embedding
$j_A: A \to \sha_C(A)^0$, is a free Baxter $C$-algebra on $A$
(of weight $\lambda$).
\end{prop}

\proof
For any $(R,P)\in \Bax_C^0$ and any morphism
$\varphi: A\to R$ in $\Alg_C^0$, we will display a unique
morphism
$\tilde{\varphi}: (\sha_C(A)^0,P_A) \to (R,P)$ in $\Bax_C^0$
that extends $\varphi$.
For each $n\in \NN$, we will define a $C$-linear map
\[ \tilde{\varphi}_n: (A^+)^{\otimes n}\otimes A \to R.\]
If $n=0$, we define 
\[ \tilde{\varphi}_0: (A^+)^{\otimes 0}\otimes A = A\to R\]
by $\tilde{\varphi}_0(x_0)=\varphi(x_0),\ x_0\in A$.
Assuming $\tilde{\varphi}_n$ is defined, we define
\[ \varphi_{n+1}:  (A^+)^{n+1}\times A \to R \]
by 
\[\varphi_{n+1} (x_0,\ldots, x_{n+1}) =
 c P (\tilde{\varphi}_n(x_1\otimes\ldots\otimes x_{n+1}))
+ \varphi(x_0\tp)
    P (\tilde{\varphi}_n(x_1\otimes\ldots\otimes x_{n+1})),\]
if $x_0=(c,x_0\tp)\in A^+=C\oplus A$. 
Using the induction hypothesis, we see that
this map is $C$-multilinear, and so induces 
\[ \tilde{\varphi}_{n+1}: (A^+)^{\otimes(n+1)}\otimes A \to R.\]
We then use $\tilde{\varphi}_n,\ n\in \NN$ to define
\begin{equation}
 \tilde{\varphi}=\sum_{n=0}^\infty \tilde{\varphi}_n:
    \sha_C(A)^0=\oplus_{n=0}^\infty
    ((A^+)^{\otimes n} \otimes A) \to R.
\label{eq:ext}
\end{equation}
Since the products of $\sha_C(A)^0\subseteq \sha_C(A^+)$
and $A$ both satisfy the mixable shuffle product identity
(see \cite[Proposition 4.2]{G-K1}), 
$\tilde{\varphi}$ is a morphism in $\Bax_C^0$.
On the other hand, for
$x_0\otimes \ldots\otimes x_n\in (A^+)^{\otimes n}\otimes A$
with $x_0=(c,x_0\tp)$, 
we have
\[x_0\otimes \ldots\otimes x_n=
 c P_A (x_1\otimes\ldots\otimes x_n)
+ \varphi(x_0\tp)
    P_A (x_1\otimes\ldots\otimes x_n). \]
Thus equation~(\ref{eq:ext}) is the only possible way to define
a morphism in $\Bax_C^0$ from $\sha_C(A)^0$ to $R$ that
extends $\varphi$.
This verifies the required universal property of $\sha_C(A)^0$.
\proofend

\section{Complete Baxter algebras}
\label{sec:comp}
We will define a natural decreasing filtration on Baxter
algebras and study the associated completion.
We will show that the completion of a shuffle Baxter algebra
is a free object in the category of complete Baxter algebras.
In this section, we retain the assumption that all algebras
have an identity. 

\subsection{Filtrations and completions}
Let $(R,P)$ be a Baxter algebra. For any subset $U$ of $R$,
denote $<U>_B$ for the Baxter ideal of $R$ generated by $U$.
We will define a decreasing filtration $\Fil^n R$
of ideals of $(R,P)$ as follows.
Define $\Fil^0 R=R$. For any $n\in \NN$, assume that
$\Fil^n R$ is defined, and inductively define
\[ \Fil^{n+1} R= <P(\Fil^n R)>_B.\]
Thus, for example,
\[ \Fil^1 R =<P(R)>_B\ {\rm\ and\ }
    \Fil^2 R=<P(<P(R)>_B)>_B. \]
Since each $\Fil^n R$ is a Baxter ideal of $R$, we have
$P(\Fil^n R)\subseteq \Fil^n R$.
Therefore,
\[\Fil^{n+1} R= <P(\Fil^n R)>_B \subseteq \Fil^n R.\]
So $\{ \Fil^n R\}_{n\in \NN}$ defines a decreasing filtration of
Baxter ideals on $(R,P)$.
Assuming $\Fil^1 R\neq R$, then each $R/\Fil^n R$, $n\in \NN_+$,
is a Baxter $C$-algebra.
Since each of the projections 
$R/\Fil^{n+1} R\to R/\Fil^n R,\ n\in \NN_+,$ is a Baxter $C$-algebra
homomorphism, the inverse limit
$\invlim (R/\Fil^n R)$ is also a Baxter $C$-algebra. 
\begin{defn}
Let $(R,P)$ be a Baxter algebra.
\begin{enumerate}
\item
The decreasing filtration $\Fil^n R$ on $R$ is called
the {\bf Baxter filtration on $R$}.
\item
The Baxter algebra $(R,P)$ is called {\bf proper} if
$\Fil^1 R$ is a proper Baxter ideal of $R$.
\item
Denote $\Bax_C\tp$ for the subcategory of $\Bax_C$ consisting
of proper Baxter $C$-algebras.
\item
For $(R,P)\in \Bax_C\tp$, 
the inverse limit $\widehat{R}=\invlim (R/\Fil^n R)$ with
the induced
Baxter operator $\widehat{P}$ is called the
{\bf (Baxter) completion of $(R,P)$}. 
\item
$(R,P)\in \Bax_C\tp$ is called {\bf (Baxter) complete} if
the natural Baxter $C$-algebra homomorphism
\[ \pi_R : R \to \invlim (R/\Fil^n R) \]
is an isomorphism.
\end{enumerate}
\end{defn}
Let $\lambda \in C$ and $A\in \Alg_C$. It is easy to see
(Proposition~\ref{prop:grade}) that 
the shuffle Baxter algebra $\sha_C(A)$ of weight $\lambda$ is proper. 
On the other hand, for $\lambda\in C$, define
$P_\lambda$ on $A$ by
$P_\lambda(a)=-\lambda a$.
Then $(A,P_\lambda)$ is a Baxter algebra of weight $\lambda$.
If $\lambda$ is invertible in $A$, then
$P_\lambda(A)=-\lambda A=A$. So $(A,P_\lambda)$ is not proper.
If $\lambda$ is not invertible in $A$, then
$(A,P_\lambda)$ is proper.
In fact, the Baxter completion of $(A,P_\lambda)$ is the same
as $\invlim A/\lambda^n A$,
the $\lambda$-adic completion of $A$.

\begin{prop}
\label{prop:comp}
For each $f: (R,P)\to (S,Q)$ in $\Bax_C\tp$, there is a unique
$\hat{f}: (\hat{R},\hat{P}) \to (\hat{S},\hat{Q})$ in $\Bax_C$ 
making the diagram
\[ \begin{array}{ccc} R &\ola{f}& S\\
    \dap{\pi_R} && \dap{\pi_S}\\
    \hat{R} &\ola{\hat{f}} & \hat{S}
    \end{array} \]
commute. 
\end{prop}
Before proving the proposition, we first give an elementary fact
on Baxter algebras.

\begin{lemma}
\label{lem:set}
Let $f: (R,P)\to (S,Q)$ be a morphism in $\Bax_C\tp$.
For any subset $U$ of $R$, we have
$f(<U>_B) \subseteq <f(U)>_B$.
\end{lemma}
\proof
For a given subset $U$ of $R$, consider the morphism
\[ R\ola{f} S\to S/<f(U)>_B\]
in $\Bax_C$. Since $U$ is in the kernel of the morphism,
it follows that $<U>_B$ is in the kernel of the morphism.
Thus $f(<U>_B)\subseteq <f(U)>_B$.
\proofend

\noindent{\bf Proof of Proposition~\ref{prop:comp}:}
We first apply induction on $k$ to show that
$f(\Fil^k R)\subset \Fil^k S$. 
For $k=0$, this just says $f(R)\subseteq S$. 
Assume that $f(\Fil^k R)\subseteq \Fil^k R$. 
Applying Lemma~\ref{lem:set} to $P(\Fil^k R)$, we obtain
\begin{eqnarray*}
f(\Fil^{k+1} R) &=& f(<P(\Fil^k R)>_B)\\
 &\subseteq& <f(P(\Fil^k R))>_B\\
& =& <Q(f(\Fil^k R))>_B\\
&\subseteq & <Q(\Fil^k S)>_B\\
&=& \Fil^{k+1} S.
\end{eqnarray*}
This completes the induction.
Then the proposition follows from general results
on completions~\cite[p.57]{Mat}. 
\proofend

\subsection{An alternative description}
We now give an interpretation of a complete free Baxter
algebra in terms infinite sequences. 
We consider the following situation.
Let $R$ be a $C$-algebra and let $R_i,\ i\in \NN$, 
be $C$-submodules of
$R$ such that
\begin{enumerate}
\item
$R= \oplus_{k\in \NN} R_k$ as a $C$-module, and 
\item
$F^n R\defeq \oplus_{k>n} R_k$ is an ideal of $R$, $n\in \NN$.
\end{enumerate}
We will define a multiplication on $\prod_{k\in \NN} R_k$.
The definition is similar to the case when $R$ is graded
$C$-algebra. For lack of a suitable reference, we give
details below.

Fix a $k\in \NN$.
Let $(x^{(n)})_n$ be a sequence of elements in $R_k$. 
If there is an $n_0\in \NN$ such that
$x^{(n)}=x^{(n_0)}$ for $n\geq n_0$,
then define
$\dislim{n\to \infty} x^{(n)} =x^{(n_0)}\in R_k$.
Further, let $(x^{(n)})=((x^{(n)}_k)_k)$ be a sequence
of elements in $\prod_{k\in \NN}R_k$.
If $\dislim{n\to \infty} x^{(n)}_k$ exists for each $k$,
then define

\[\dislim{n\to \infty} x^{(n)}
    = (\dislim{n\to \infty} x^{(n)}_k)_k \in \prod_k R_k.\]
For any $x=(x_k)_k\in \prod_k R_k$ and any $n\in \NN$, define
$x^{[n]} =(x^{[n]}_k)_k \in \prod_{k\in\NN} R_k$ by 
\[ x^{[n]}_k = \left \{ \begin{array}{ll}
    x_k, & k\leq n,\\
    0, & k>n \end{array} \right .\]
Now let $x=(x_k)_k$ and $y=(y_k)_k$ be two elements of
$\prod_{k\in \NN} R_k$.
For given $n\in \NN$, we have $x^{[n]},\ y^{[n]}\in
    R=\oplus_k R_k$.
So $x^{[n]} y^{[n]}$ can be uniquely expressed as
$(z^{(n)}_k)_k$, $z^{(n)}_k\in R_k$ and
$z^{(n)}_k=0$ for $k>>0$.
For each fixed $k$, we obtain a sequence
$(z^{(n)}_k)_n$ in $R_k$.
When $n\geq k$, we have
$x^{[n]}=x^{[k]}+x\tp$
and
$y^{[n]}=y^{[k]}+y\tp$
with $x\tp, y\tp \in F^{k+1} R$.
Since $F^{k+1} R$ is an ideal of $R$, we have
\[ x^{[n]} y^{[n]}\equiv x^{[k]}y^{[k]}  \pmod{F^{k+1}R}.\]
Therefore
$z^{(n)}_k =z^{(k)}_k$ for $n\geq k$ and
$\dislim{n\to \infty} z^{(n)}_k \in R_k$ is well-defined.
Define
\[ (x_k)_k (y_k)_k =(\dislim{n\to \infty} z^{(n)}_k)_k
    \in \prod_k R_k.\]
In other words,
\[ (x_k)_k (y_k)_k = \dislim{n\to \infty} x^{[n]}y^{[n]}.\]
It can be easily verified that this defines an associative,
commutative
multiplication on $\prod_{k\in \NN}R_k$, making it into a
commutative $C$-algebra.

\begin{prop}
\label{prop:series}
Let $R$ be a $C$-algebra and let $R_i,\ i\in \NN$, 
be $C$-submodules of
$R$ such that
\begin{enumerate}
\item
$R= \oplus_{k\in \NN} R_k$ as a $C$-module, and 
\item
$F^n R\defeq \oplus_{k>n} R_k$ is an ideal of $R$, $n\in \NN$.
\end{enumerate}
There is a unique $C$-algebra isomorphism 
\[ \psi_R:
 \invlim (R/F^k R) \cong \prod_{k\in \NN} R_k\]
that makes the diagram in $\Alg_C$

\[ \begin{array}{ccc}
R & \hookrightarrow & \prod_{k\in \NN} R_k \\
\dap{\pi_R} & \ \ \  \nearrow_{\psi_R} \!\!\! & \\
\invlim (R/F^k R) && \end{array} \]
commute. 
\end{prop}

\proof
By definition, $\invlim (R/F^k R)$ is the inverse limit of the
inverse system $p_{n+1,n}:R/F^{n+1} R\to R/F^n R,\ n\geq 1$. 
An element
$((x^{(n)}_k)_k +F^n R )_n\in \prod_{n\in \NN} R/F^n R$
is an element of $\invlim (R/F^k R)$ if and only if,
for any $n\geq 1$, 
\[ p_{n+1,n}((x^{(n+1)}_k)_k +F^{n+1} R)=
    (x^{(n)}_k)_k +F^n R.\]
Since $R/F^n R\cong \oplus_{k\leq n}R_k$, this is so if and only
if $x^{(n+1)}_k=x^{(n)}_k$ for $k\leq n$.
Therefore $((x^{(n)}_k)_k + F^n R)_n$ is in $\invlim (R/F^k R)$
if and only if there is $(y_n)_n\in \prod_{n\in \NN} R_n$
such that, for any $n$,
$x^{(n)}_k = y_k$ for $k\leq n$.
In fact, we can take $y_k=x^{(k)}_k$. 
This gives the desired map
$\psi_R: \invlim (R/F^k R)\to \prod_{k\in \NN}R_k$.
More precisely, we have 

\begin{equation}
 \psi_R(((x^{(n)}_k)_k + F^n R)_n)
= (x^{(k)}_k)_k.
\label{eq:series}
\end{equation}
For $x=(x_k)_k \in R=\oplus_{k\in \NN} R_k$,
$\pi_R(x)\in \invlim (R/F^k R)$ is defined to be the sequence
$((x_k)_k + F^n R )_n$ which corresponds under $\psi_R$ to the
element $(x_k)_k\in \prod_{k\in \NN}R_k$.
This proves the commutativity of the diagram.
\proofend

\subsection{Complete shuffle Baxter algebras}
We now consider the completion of $\sha_C(A)$. 
Recall that we denote $\sha^k_C(A)$ for the $C$-submodule
$A^{\otimes (k+1)}$ of $\sha_C(A)$.
We denote $\invlim \sha_C(A)/\Fil^k \sha_C(A)$
by $\widehat{\sha}_C(A)$. 

\begin{prop}
\label{prop:grade}
Given $k\in \NN_+$, 
\begin{enumerate}
\item
$\Fil^k \sha_C(A)= \bigoplus_{n\geq k} \sha^n_C(A),$ 
\item
$\Fil^k \sha_C(A)$ is a Baxter homogeneous
ideal of $\sha_C(A)$, and 
\item
The quotient Baxter $C$-algebra $\sha_C(A)/\Fil^k\sha_C(A)$ 
is isomorphic to $\bigoplus_{n=0}^{k-1} \sha^n_C(A)$ as
a $C$-module.
\end{enumerate}
\end{prop}

\proof
1.
It follows from the definition of $\shpr$ that,
for any $C$-algebra $A$, 
\[ \sha^m_C(A) \shpr  \sha^n_C(A) \subseteq
    \sum_{k=\max\{m,n\}}^{m+n}\sha^k_C(A).\]
This shows that 
$ \bigoplus_{n\geq k} \sha_C^n(A)$ is a Baxter ideal. 
Next we prove

\begin{equation}
\Fil^k \sha_C(A)= \bigoplus_{n\geq k} \sha^n_C(A)
\label{eq:fil}
\end{equation}
by induction on $k$. 
By definition, 
$\Fil^1 \sha_C(A)$ is the Baxter ideal generated by
$P_A(\sha_C(A))=\bfone_A \otimes \sha_C(A)$. 
On the other hand,
$\bigoplus_{n\geq 1} \sha_C^n(A)$ equals 
$A\shpr (\bfone_A \otimes \sha_C(A))$,
hence is also generated by $\bfone_A \otimes \sha_C(A)$.
This verifies equation~(\ref{eq:fil}) for $k=1$. 
Assume that
$\Fil^k \sha_C(A)= \bigoplus_{n\geq k} \sha^n_C(A)$
for a $k\in \NN^+$.  From this we obtain that 
$\Fil^{k+1} \sha_C(A)$ is the Baxter ideal generated by
\[ P_A(\Fil^k \sha_C(A))=P_A(\bigoplus_{n\geq k} \sha^n_C(A))
= \bfone_A\otimes (\bigoplus_{n\geq k} \sha^n_C(A)).\]
On the other hand,
$\bigoplus_{n\geq k+1} \sha_C^n(A)$ equals 
$R\shpr (\bfone_A \otimes \bigoplus_{n\geq k} \sha_C^n(A))$,
hence is also generated by
$\bfone_A \otimes \bigoplus_{n\geq k} \sha_C^n(A)$.
This verifies equation~(\ref{eq:fil}) for $k+1$.

Other statements in the proposition follow immediately from
the first one.
\proofend

It follows from Proposition~\ref{prop:grade}
that $R=\sha_C(A)=\oplus_{k\in \NN} \sha_C^k (A)$ 
satisfies the two conditions for $R$ in
Proposition~\ref{prop:series}. 
Thus the product
$\prod_{k\in \NN} \sha_C^k (A)$ is a $C$-algebra.
Define an operator $\bar{P}$ on this product by
$\bar{P}((x_k)) = (\bfone_A\otimes x_{k-1})$
with the convention that $\bfone_A\otimes x_{k-1}=0$ for $k=0$.

\begin{theorem}
\label{thm:series}
\begin{enumerate}
\item
$\bar{P}$ is a Baxter operator on
$\prod_{k\in \NN} \sha_C^k (A)$, and 
$\psi_{A}\defeq \psi_{\sha_C(A)}:
\widehat{\sha}_C(A) \to \prod_{k\in \NN} \sha_C^k (A)$
is an isomorphism of $C$-Baxter algebras.
\item
Given a morphism $f:A\to B$ in $\Alg_C$, we have the following
commutative diagram in $\Bax_C$ 
\[ \begin{array}{ccc}
\widehat{\sha}_C(A) & \ola{\psi_A} &
    \prod_{k\in\NN} \sha_C^k (A) \\
\dap{\widehat{\sha}_C(f)} &&    \dap{\prod_k f_k} \\
\widehat{\sha}_C(B) & \ola{\psi_B} &
    \prod_{k\in\NN} \sha_C^k (B)
\end{array} \]
where
$\widehat{\sha}_C(f)$ is induced by $\sha_C(f)$ which is
in turn induced by $f$,
and 
$f_k:\sha_C^k(A) \to \sha_C^k(B)$ is the
tensor power morphism of $C$-modules
$f^{\otimes (k+1)}: A^{\otimes (k+1)} \to B^{\otimes (k+1)}$
induced from $f$. 
\end{enumerate}
\end{theorem}

\proof
1.
Let 
$((x^{(n)}_k)_k + \Fil^n \sha_C(A))_n\in \widehat{\sha}_C(A)$
be given. 
Using the formula~(\ref{eq:series}) for the map $\psi_A$, we have

\begin{eqnarray*}
\lefteqn{
(\bar{P}\circ \psi_A)
(((x^{(n)}_k)_k + \Fil^n \sha_C(A) )_n) }\\
&=& \bar{P}((x^{(k)}_k)_k) \\
&=& (\bfone_A \otimes x^{(k)}_{k-1})_k, 
\end{eqnarray*}
and

\begin{eqnarray*}
\lefteqn{(\psi_A\circ \hat{P})
(((x^{(n)}_k)_k + \Fil^n \sha_C(A) )_n) }\\
&=&
\psi_A( ((\bfone_A\otimes x^{(n)}_{k-1})_k
    + \Fil^n \sha_C(A) )_n)\\
&=&(\bfone_A\otimes x^{(k)}_{k-1})_k\\
\end{eqnarray*}
Thus
$\bar{P}\circ \psi_A = \psi_A \circ \hat{P}$. 
Since $\psi_A$ is an isomorphism of $C$-algebras and
since $\hat{P}$ is known to satisfy the identity defining 
a Baxter operator, the above equation implies that the same
identity is satisfied by $\bar{P}$. 

\noindent
2.
Given
$((x^{(n)}_k)_k + \Fil^n \sha_C(A) )_n\in \widehat{\sha}_C(A)$,
using the formula~(\ref{eq:series}) we have

\begin{eqnarray*}
\lefteqn{
(\prod_k f_k\circ \psi_A)
(((x^{(n)}_k)_k + \Fil^n \sha_C(A) )_n) }\\
&=&
\prod_k f_k (( x^{(k)}_k)_k) \\
&=& (f_k( x^{(k)}_k))_k
\end{eqnarray*}
and

\begin{eqnarray*}
\lefteqn{
(\psi_A\circ \widehat{\sha}_C(A))
(((x^{(n)}_k)_k + \Fil^n \sha_C(A) )_n) }\\
&=&
\psi_A(((f_k(x^{(n)}_k))_k + \Fil^n \sha_C(B) )_n) \\
&=&
( f_k(x^{(k)}_k))_k. 
\end{eqnarray*}
This proves that the diagram commutes. 
\proofend

As an example of Theorem~\ref{thm:series},
consider the case when $A=C$ and $\lambda=0$.
Let $HC$ be the ring of Hurwitz series over $C$~\cite{Ke},
defined to be the set of sequences 
\[ \{ (a_n) \mid a_n\in C, n\in \NN \} \]
in which the addition
is defined componentwise and the multiplication is
defined by
\[ (a_n) (b_n) =(c_n)\]
with
\[ c_n =\displaystyle{ \sum_{k=0}^n }
    \bincc{n}{k}a_kb_{n-k}.\]
Denote $e_n$ for the sequence $(a_k)$ in which
$a_n=\bfone_C$ and $a_k=0$ for $k\neq n$.
Since $e_n e_m = \bincc{m+n}{n}e_{m+n}$,
the following corollary follows from Proposition~\ref{prop:unit}. 

\begin{coro}
\label{co:hurw}
The assignment
\[ \bfone^{\otimes (n+1)} \mapsto e_n,\ n\geq 0\]
defines an isomorphism
\[ \widehat{\sha}_C(C) \to HC.  \]
\end{coro}
\proofend

By Theorem~\ref{thm:series} and part three of
Proposition~\ref{prop:grade},
we have the isomorphism of inverse systems
\[ \widehat{\sha}_C(A)/ \Fil^k \widehat{\sha}_C(A)
    \cong \sha_C(A)/ \Fil^k \sha_C(A). \]
Thus the completion of $\widehat{\sha}_C(A)$ is itself,
so it is complete. 
We next verify the free universal property of
$\widehat{\sha}_C(A)$ in the category $\widehat{\Bax}_C$
of complete Baxter $C$-algebras. 
Abbreviate $\pi_A$ for $\pi_{\sha_C(A)}$.

\begin{theorem}
\label{thm:comp}
$(\widehat{\sha}_C(A),\hat{P}_A)$, together with the natural embedding
$\hat{j}_A:A\ola{j_A} \sha_C(A)\ola{\pi_A}
    \widehat{\sha}_C(A)$,
is a free complete Baxter $C$-algebra on $A$ (of weight $\lambda$). 
In other words,
for any complete Baxter $C$-algebra $(R,P)$ and any
$C$-algebra map
$\varphi:A\rar R$, there exists
a unique Baxter $C$-algebra homomorphism
$\hat{\varphi}:(\widehat{\sha}_C(A), \hat{P}_A)\rar (R,P)$
such that the  diagram
\[\xymatrix{
A \ar[rr]^(0.4){\hat{j}_A} \ar[drr]_{\varphi}
    && \widehat{\sha}_C(A) \ar[d]^{\hat{\varphi}} \\
&& R } \]
commutes. 
\end{theorem}

\proof
Given a $(R,P)$ and $\varphi: A\to R$ as in the statement
of the theorem, by the universal property of
$(\sha_C(A),j_A)$ in $\Bax_C$, there is a unique
$\tilde{\varphi}:\sha_C(A) \to R$ in $\Bax_C$ such that
$\tilde{\varphi}\circ j_A=\varphi$.
Since $R$ is complete, by Proposition~\ref{prop:comp}, there is
\[ \hat{\varphi}\defeq \hat{\tilde{\varphi}}:
\widehat{\sha}_C(A) \to \hat{R}\cong R \]
such that
$\hat{\varphi}\circ \pi_A = \tilde{\varphi}$.
Then we have

\[
\hat{\varphi}\circ \widehat{j}_A 
= \hat{\varphi}\circ \pi_A \circ j_A
= \tilde{\varphi}\circ j_A 
= \varphi. \]

This proves the existence of $\hat{\varphi}$.
The uniqueness follows from the uniqueness of
$\tilde{\varphi}$ and the uniqueness of
the completion.
\proofend

\subsection{Completions of Baxter algebras}
It is clear that complete Baxter $C$-algebras,
together with the Baxter algebra homomorphisms between them,
form a full subcategory $\widehat{\Bax}_C$ of $\Bax_C$.
Recall that we denote $\Bax_C\tp$ for the full subcategory
of $\Bax_C$ consisting of proper Baxter algebras. 
Denoting $I_C: \widehat{\Bax}_C \to \Bax_C\tp$
for the natural inclusion of categories, 
we then have 

\begin{prop}
\label{prop:comp2}
\begin{enumerate}
\item
For any $R\in \Bax_C\tp$,
the Baxter completion $\hat{R}$ of $R$ is complete.
\item
The assignments $(R,P)\mapsto (\hat{R},\hat{P})$ and 
$f \mapsto \hat{f}$ 
define a functor $F_C$ 
from $\Bax_C\tp$ to $\widehat{\Bax}_C$, and the 
morphisms $\pi_R:R\to \hat{R}$, $R\in \Bax_C\tp$, 
define a natural transformation
between the identity functor on $\Bax_C\tp$ and the
functor $I_C\circ F_C:\Bax_C\tp \to \Bax_C\tp$.
\end{enumerate}
\end{prop}
\proof
We only need to prove that $\hat{R}$ is complete. 
The rest of the proof is clear.

Note that the filtration $\Fil^k$ on $\hat{R}$ is not defined to be
the natural filtration induced by the filtration $\Fil^k$ on $R$.
So the general results on completions do not apply.
Instead, we will use the facts that any Baxter algebra is a quotient
of a shuffle Baxter algebra, and that the completion of a shuffle Baxter algebra
is complete, which follows from Proposition~\ref{prop:grade} and
Theorem~\ref{thm:series}. 

For a given $(R,P)\in \Bax_C\tp$, by the universal property of
free Baxter algebras, there is an $A\in \Alg_C$ and a surjective
morphism $\psi: \sha_C(A)\to R$ in $\Bax_C$.
We will prove by induction on $n\in \NN^+$ that

\begin{equation}
 \psi(\Fil^n \sha_C(A)) = \Fil^n R.
\label{eq:fil2}
\end{equation}
For $n=1$ we have
\begin{eqnarray*}
\lefteqn{ \psi(\Fil^1 \sha_C(A))= \psi(\sha_C(A)P_A(\sha_C(A)))}\\
&=& \psi(\sha_C(A))P(\psi(\sha_C(A)))\\
&=& R P(R).
\end{eqnarray*}
Note that $R P(R)$ is the ideal of $R$ generated by $P(R)$.
Since
$P(R P(R))\subseteq P(R)$, it is in fact the Baxter ideal of $R$
generated by $P(R)$. Thus $R P(R)=\Fil^1 R$.
So equation~(\ref{eq:fil2}) holds for $n=1$.

Assume that the equation holds for $n$.
Part one of Proposition~\ref{prop:grade}
shows that 
$\Fil^{n+1}\sha_C(A)$ is the ideal of $\sha_C(A)$ generated by
$P_A(\Fil^n \sha_C(A))$.
Then by induction we have
\begin{eqnarray*}
\lefteqn{ \psi(\Fil^{n+1} \sha_C(A))
    =\psi (\sha_C(A) P_A(\Fil^n \sha_C(A)))}\\
&=& \psi(\sha_C(A))P(\psi(\sha_C(A)))\\
&=& R P(\Fil^n R).
\end{eqnarray*}
Since
\begin{eqnarray*}
\lefteqn{ P(R P(\Fil^n R))
    =P(\psi(\sha_C(A))P(\psi(\Fil^n \sha_C(A))))}\\
&=& \psi(P_A(\sha_C(A)P_A(\Fil^n \sha_C(A))))\\
&=& \psi(P_A(\Fil^{n+1}\sha_C(A)))\\
&\subseteq & \psi(\sha_C(A)P_A(\Fil^n \sha_C(A)))\\
&=& RP(\Fil^n R),
\end{eqnarray*}
$R P(\Fil^n R)$ is the Baxter ideal of $R$ generated by
$P(\Fil^n R)$, so is equal to $\Fil^{n+1} R$.
This completes the induction.

Because of equation~(\ref{eq:fil2}), the morphism
$\psi: \sha_C(A) \to R$ induces a morphism
\[ \psi_n: \sha_C(A)/\Fil^n \sha_C(A) \to R/\Fil^n R\]
for each $n\in\NN_+$ and the kernel of $\psi_n$ is
$(\ker \psi +\Fil^n \sha_C(A))/\Fil^n \sha_C(A)$ which is isomorphic to
$\ker \psi/(\ker \phi \cap \Fil^n \sha_C(A))$.
Thus we have the exact sequence of inverse systems
\[ 0\!\! \to (\ker \psi +\Fil^n \sha_C(A))/\Fil^n \sha_C(A)
\to \sha_C(A)/ \Fil^n \sha_C(A)
\to R/\Fil^n R \to\!\! 0\]
and the transition map of the left inverse system is identified with
the natural map
\[ \ker \psi/(\ker \psi \cap \Fil^{n+1} \sha_C(A))
\to \ker \psi/(\ker \psi \cap \Fil^n \sha_C(A))\]
so is surjective.
By~\cite[Lemma 3.5.3]{We}, for the first derived functor
$R^1 \invlim $ of the inverse limit, 
\[ R^1 \invlim (\ker \psi +\Fil^n \sha_C(A))/\Fil^n \sha_C(A) =0.\]
Therefore the above exact sequence of inverse systems gives
the surjective morphism
\[\hat{\psi}: \widehat{\sha}_C(A) \to \hat{R}.\]

Because of Theorem~\ref{thm:series}, 
$ \Fil^{n+1} \widehat{\sha}_C(A)$ is the ideal of
$\widehat{\sha}_C(A)$ generated by
$\hat{P}_A(\Fil^n \widehat{\sha}_C(A))$. 
Then the same argument for
$\psi: \sha_C(A) \twoheadrightarrow R$
in the previous part of the proof
can be repeated for the morphism
$\hat{\psi}: \widehat{\sha}_C(A) \twoheadrightarrow \hat{R}$. 
In particular, we have, for any $n\in \NN^+$,

\begin{equation}
 \psi(\Fil^n \widehat{\sha}_C(A)) = \Fil^n \hat{R}.
\label{eq:fil3}
\end{equation}
We then obtain a surjective
morphism 
\[ \hat{\hat{\psi}}: \widehat{\widehat{\sha}}_C(A) \to \hat{\hat{R}}.\]
Since $\widehat{\sha}_C(A)$ is its own completion, we have the
commutative diagram
\[ \begin{array}{ccc}
\widehat{\sha}_C(A) & \cong & \widehat{\widehat{\sha}}_C(A)\\
\dap{\hat{\psi}} && \dap{\hat{\hat{\psi}}}\\
\hat{R} & \ola{\pi_{\hat{R}}} &\hat{\hat{R}}
\end{array}\]
in which both of the vertical maps are surjective.
By the commutativity of the diagram, $\pi_{\hat{R}}$ is
surjective. 
Since $\cap_n \Fil^n \widehat{\sha}_C(A) =0$, 
by equation~(\ref{eq:fil3}), we have 
$\cap_n \Fil^n \hat{R}=0$.
Thus $\pi_{\hat{R}}$ is injective.
Therefore, $\hat{R}$ is complete. 
\proofend

\section{The standard Baxter algebra}
\label{sec:rota}
The standard Baxter algebra constructed by Rota in~\cite{Rot1}
is a free object in
the category $\Bax_C^0$ of Baxter algebras not necessarily
having an identity. It is described as a Baxter subalgebra
of another Baxter algebra whose construction goes back to
Baxter~\cite{Ba}.
In Rota's construction, there are further restrictions 
that $C$ is a field of characteristic zero,
the free Baxter algebra obtained is
on a finite set $X$, and the weight $\lambda$ is $1$.
By making use of shuffle Baxter algebras, 
we will show that Rota's description can be modified to yield 
a free Baxter algebra on an algebra in the category $\Bax_C$
of Baxter algebras with an identity,
with a mild restriction on the weight $\lambda$. 
We will also provide a similar construction for
algebras not necessarily having an identity,
and for complete Baxter algebras. 

\subsection{The standard Baxter algebra of Rota}
We first briefly recall the construction of Rota of
a standard Baxter algebra
$\frakS(X)$ on a set $X$. For details, see~\cite{Rot1,R-S}.

As before, let $C$ be a commutative ring with an identity, 
and fix $\lambda \in C$. Let $X$ be a given set.
For each $x\in X$, let $t^{(x)}$ be a sequence
$t^{(x)}=( t^{(x)}_1,\ldots, t^{(x)}_n,\ldots )$
of distinct symbols $t^{(x)}_n$.
We also require that the sets $\{t^{(x_1)}_n\}_n$ and
$\{t^{(x_2)}_n\}_n$ are disjoint for $x_1\neq x_2$ in $X$.
Denote
\[ \overline{X} = \cup_{x\in X} \{t^{(x)}_n \mid n\in \NN_+\} \]
and denote $\frakA (X)$ for the ring of sequences with entries in 
$C[\overline{X}]$, the $C$-algebra of polynomials with variables in
$\overline{X}$.
Thus the addition, multiplication and scalar multiplication by
$C[\overline{X}]$ in $\frakA(X)$ are defined componentwise.
It will useful to have the following description of $\frakA(X)$. 
For $k\in \NN_+$, denote $\gamma_k$ for the sequence
$(\delta_{n,k})_n$, where $\delta_{n,k}$ is the Kronecker delta. 
Then we can identify a sequence $(a_n)_n$ in $\frakA(X)$
with a series
\[ \sum_{n=1}^\infty a_n \gamma_n = a_1 \gamma_1 +a_2 \gamma_2 + \ldots .\]
Then the addition, multiplication and scalar multiplication by
$C[\overline{X}]$ are given termwise. 

Define
\[P_X\tp=P_{X,\lambda}\tp: \frakA(X)\to \frakA(X)\]
by
\[ P_X\tp(a_1,a_2,a_3,\ldots)
=\lambda (0,a_1,a_1+a_2,a_1+a_2+a_3,\ldots).\]
In other words, each entry of $P_X\tp (a),\ a=(a_1,a_2,\ldots),$
is $\lambda$ times the sum of the previous entries of $a$.
If elements in $\frakA(X)$ are described by series
$\sum_{n=1}^\infty a_n \gamma_n$ given above, then we simply have
\[ P_X\tp(\sum_{n=1}^\infty a_n \gamma_n)=
    \lambda \sum_{n=1}^\infty (\sum_{i=1}^{n-1} a_i) \gamma_n.\]
It is well-known~\cite{Ba,Rot1} that, for $\lambda=1$,
$P_X\tp$ defines a Baxter operator of weight $1$ on $\frakA(X)$.
It follows that, for any $\lambda\in C$,
$P_X\tp$ defines a Baxter operator of weight $\lambda$ on $\frakA(X)$,
since it can be easily verified that 
for any Baxter operator $P$ of weight $1$,  the operator
$\lambda P$ is a Baxter operator of weight $\lambda$.
Hence $(\frakA(X),P\tp_X)$ is in $\Bax_C$. 

\begin{defn}
Let $\frakS(X)^0$ be the Baxter subalgebra in $\Bax_C^0$
of $\frakA(X)$ generated
by the sequences $t^{(x)}=(t^{(x)}_1,\ldots,x^{(x)}_n,\ldots),\ x\in X$.
$\frakS(X)^0$ is called the {\bf standard Baxter algebra}
on $X$. 
\end{defn}
Note that 
$\frakS(X)^0$ is denoted by $\frakS(X)$ in Rota's notation.
We reserve $\frakS(X)$ for the free Baxter algebra on $X$ with
an identity that will be defined below.

\begin{theorem}
\label{thm:rota}
{\bf (Rota)}~\cite{Rot1}
$(\frakS(X)^0,P_X\tp)$ is a free Baxter algebra on $X$ in the
category $\Bax_C^0$.
\end{theorem}

\subsection{The standard Baxter algebra in general}
Given $A\in \Alg_C$, 
we now give an alternative construction of
a free Baxter algebra on $A$ in the category $\Bax_C$
of Baxter algebras with an identity.

For each $n\in \NN_+$, denote $A^{\otimes n}$ for the tensor power
algebra.
Denote the direct limit algebra 
\[ \overline{A}= \dirlim A^{\otimes n} \]
where the transition map is given by
\[ A^{\otimes n} \rightarrowtail A^{\otimes (n+1)},\
    x \mapsto x\otimes \bfone_A.\]
Note that the multiplication on $A^{\otimes n}$ here is
different from the multiplication on $A^{\otimes n}$ when it
is regarded as the $C$-submodule $\sha_C^{n-1}(A)$ of $\sha_C(A)$.
To distinguish between the two contexts, we will use the notation
$\sha_C^n(A)$ for $A^{\otimes (n+1)}\subseteq \sha_C(A)$. 
Let $\frakA(A)$ be the set of sequences with entries in
$\overline{A}$. Thus we have
\[ \frakA(A) = \prod_{n=1}^\infty \overline{A} \gamma_n
    =\left \{ \sum_{n=1}^\infty a_n \gamma_n,
    a_n \in \overline{A} \right \}.\]
Define addition, multiplication and scalar multiplication
on $\frakA(A)$ componentwise, making $\frakA(A)$ into a
$\overline{A}$-algebra,
with the sequence $(1,1,\ldots)$ as the identity. 
Define
\[P_A\tp=P_{A,\lambda}\tp: \frakA(A)\to \frakA(A)\]
by
\[ P_A\tp(a_1,a_2,a_3,\ldots)
=\lambda (0,a_1,a_1+a_2,a_1+a_2+a_3,\ldots).\]
Then $(\frakA(A),P\tp_A)$ is in $\Bax_C$. 
For each $a\in A$, define $t^{(a)}=(t^{(a)}_k)_k$ in $\frakA(A)$ by 
\[ t^{(a)}_k=\otimes_{i=1}^k a_i
    (=\otimes_{i=1}^\infty a_i),\ a_i = \left \{
    \begin{array}{ll} a, & i=k,\\
    1, & i\neq k. \end{array} \right .  \]

\begin{defn}
Let $\frakS(A)$ be the Baxter subalgebra in $\Bax_C$
of $\frakA(A)$ generated
by the sequences
$t^{(a)}=(t^{(a)}_1,\ldots,t^{(a)}_n,\ldots),\ a\in A$.
$\frakS(A)$ is called the {\bf standard Baxter algebra}
on $A$. 
\end{defn}

Since $\sha_C(A)$ is a free Baxter algebra on $A$,
the morphism in $\Alg_C$ 
\[ A \to \frakA(A), a\mapsto t^{(a)}\]
extends uniquely to a morphism in $\Bax_C$

\[ \Phi: \sha_C(A) \to \frakA(A). \]
We will prove in Theorem~\ref{thm:s-r} that, when $\lambda$ 
is not a zero divisor in $\overline{A}$, $\Phi$ is an isomorphism.
Hence $(\frakS(A),P\tp_A)$ is a free Baxter algebra
on $A$ in the category $\Bax_C$. 
Before proving the theorem, 
we will first give some notations and preliminary results. 

For $k\in \NN_+$, denote 
\[F^k \frakA (A) = \{ (a_i)\in \frakA(A)\mid a_i=0,\ i\leq k\}
=\{  \sum_{n=k+1}^\infty a_n \gamma_n \mid a_n\in \overline{A} \}.\]
Also denote $F^0 \frakA(A)=\frakA(A)$. 
Clearly each $F^k \frakA(A)$ is an ideal of $\frakA(A)$.
Define 
\[ F^k \frakS(A)=F^k \frakA(A)\cap \frakS(A). \]
Then we have
\[ F^k \frakS(A) = \{ (a_i)\in \frakS(A)\mid a_i=0,\ i\leq k\}.\]
$F^k \frakS(A)$ are ideals of $\frakS(A)$. 
Recall from Section~\ref{sec:comp}
that there is a canonical (Baxter) filtration $\Fil^k$ on
$\frakA(A)$ and $\frakS(A)$ defined by $P\tp_A$. 
We will explain the relation between them in
Lemma~\ref{lem:fil2}. 

\begin{lemma}
\label{lem:fil}
For any $k\in \NN$, we have
\begin{enumerate}
\item
$\frakA(A) P_A\tp(F^k \frakA(A)) \subseteq F^{k+1} \frakA(A)$. 
\item
$\Phi(\Fil^k \sha_C(A)) \subseteq F^k \frakA(A)$.
\end{enumerate}
Similar inclusions hold for $\frakS(A)$. 
\end{lemma}
\proof
We only need to verify the inclusions for $\frakA(A)$.
The inclusions for $\frakS(A)$ follows immediately.
By the definition of $P_A\tp(F^k \frakA(A))$ we have 
$P_A\tp(F^k \frakA(A)) \subseteq F^{k+1} \frakA(A)$.
Since $\frakA(A) P_A\tp(F^k \frakA(A))$ equals to the
ideal of $\frakA(A)$ generated by $P_A\tp(F^k \frakA(A))$,
we get the first inclusion. 

The second inclusion is clear for $k=0$.
By induction on $k$, for $k>0$, we have 

\begin{eqnarray*}
\lefteqn{ \Phi(\Fil^k \sha_C(A)) \subseteq
\Phi (\sha_C(A) P_A(\Fil^{k-1}\sha_C(A)))}\\
&\subseteq &
\Phi (\sha_C(A)) P_A\tp (\Phi(\Fil^{k-1}\sha_C(A))) \\
&\subseteq &
\frakA(A) P_A\tp(F^{k-1} \frakA(A)) \\
&\subseteq &
F^k \frakA(A). 
\end{eqnarray*}
\proofend

\begin{lemma}
\label{lem:lead}
For $n\in \NN_+$ and $a_1\otimes\ldots \otimes a_n\in \sha_C^n(A)$, 
we have 
\[\Phi (a_1\otimes \ldots \otimes a_n)
    =\lambda^{n-1}(a_n \otimes \ldots \otimes a_1) \gamma_n
    + F^n \frakS(A).\]
\end{lemma}

\proof
By definition, for $a_1\in A\subseteq \sha_C(A)$, 
\[ \Phi (a_1)= \sum_{k=1}^\infty t_k^{(a_1)} \gamma_k
    = a_1\gamma_1+(1\otimes a_1)\gamma_2 + \ldots . \]
So the lemma is proved for $n=1$.
Assume that the lemma is proved for $n$, and consider
$a_1\otimes \ldots \otimes a_{n+1}\in \sha_C^n(A)$. 
Applying Lemma~\ref{lem:fil}, we have
\begin{eqnarray*}
\lefteqn{\Phi(a_1\otimes a_2 \otimes 
    \ldots \otimes a_{n+1})
= \Phi(a_1 P_A(a_2\otimes\ldots \otimes a_{n+1})) }\\
&=& \Phi(a_1)\Phi (P_A(a_{n+1}\otimes\ldots\otimes a_2))\\
&=& \Phi(a_1)P\tp_A(\Phi(a_{n+1}\otimes\ldots\otimes a_2))\\
&=& ( \sum_{k=1}^\infty t^{(a_1)}_k \gamma_k  )\ 
P_A\tp (\lambda^{n-1}(a_{n+1}\otimes \ldots \otimes a_2) \gamma_n
    +{\rm\ a\ term\ in\ } F^n \frakS(A))\\
&=& ( \sum_{k=1}^\infty t^{(a_1)}_k \gamma_k )
 (\lambda^n (a_{n+1}\otimes \ldots \otimes a_2) \gamma_{n+1}
    +{\rm\ a\ term\ in\ } F^{n+1}\frakS(A))\\
&=& \lambda^n(a_{n+1}\otimes a_n\otimes \ldots \otimes a_1)\gamma_{n+1}
    +{\rm\ a\ term\ in\ } F^{n+1}\frakS(A).
\end{eqnarray*}
This completes the induction.
\proofend

Now we are ready to prove the main theorem of this section. 

\begin{theorem}
\label{thm:s-r}
Assume that $\lambda\in C$ is not a zero divisor in $\overline{A}$. 
The morphism in $\Bax_C$ 
\[ \Phi: \sha_C(A)\to \frakS(A) \]
induced by sending
$a\in A$ to $t^{(a)}=(t^{(a)}_1,\ldots,t^{(a)}_n,\ldots)$
is an isomorphism.
\end{theorem}

\begin{coro}
When $\lambda$ is not a zero divisor in $\overline{A}$, 
$(\frakS(A),P_A\tp)$ is a free Baxter algebra on $A$ in the
category $\Bax_C$. \proofend
\end{coro}

\begin{coro}
\label{co:rx}
Assume that $\lambda$ is not a zero divisor in $C$. 
Let $X$ be a set. 
The morphism in $\Bax_C$ 
\[ \Phi: \sha_C(X)\to \frakS(X) \]
induced by sending
$x\in X$ to $t^{(x)}=(t^{(x)}_1,\ldots,t^{(x)}_n,\ldots)$
is an isomorphism.
The restriction of $\Phi$ to
$\sha_C(X)^0$ is an isomorphism in $\Bax_C^0$
from $\sha_C(X)^0$ to $\frakS(X)^0$. 
\end{coro}
\proof
Applying Theorem~\ref{thm:s-r} to the case when $A=C[X]$,
we obtain
$\Phi: \sha_C(X)\cong \frakS(X)$. 
Since $\sha_C(X)^0\subseteq \sha_C(X)$ is generated by $X$
in $\Bax_C^0$ and
$\frakS(X)^0\subseteq \frakS(X)= \Phi(\sha_C(X))$
is generated by $\Phi(X)$ in $\Bax_C^0$, the corollary follows.
\proofend

\remarks

\noindent
1. 
The proof of Theorem~\ref{thm:s-r}, specialized to the setting
of Corollary~\ref{co:rx}, also gives another proof of
Theorem~\ref{thm:rota}.

\noindent
2. 
The above construction of $\frakS(A)$ for $A\in \Alg_C$
can be modified to give the construction of an internal
free Baxter algebra $\frakS(A)^0$ in $\Bax_C^0$ on $A$
for $A\in \Alg_C^0$.
The situation is similar to the construction of
shuffle Baxter algebras not necessarily having an identity, 
discussed in section~\ref{sec:back}.

\vspace{.2cm}

\noindent
{\bf Proof of Theorem~\ref{thm:s-r}:}
Since $(\sha_C(A),P_A)$ is a free Baxter algebra on $A$
in $\Bax_C$, the assignment
\[ \Phi: A \to \frakS(A),\
a\mapsto t^{(a)}=(t^{(a)}_1,\ldots,t^{(a)}_n,\ldots)\]
induces a morphism
$\Phi: \sha_C(A)\to \frakS(A)$ in $\Bax_C$.
Since $\frakS(A)$ is the Baxter subalgebra of $\frakA(A)$
generated by $A$, the morphism $\Phi$ is
onto.
So we only need to verify that $\Phi$ is injective.

For any $G\in \sha_C(X)$, we can uniquely write
$G=\sum_{n\in \NN} G_n$ with $G_n\in \sha_C^n(A)$.
Suppose $\Phi(G)=0$; we will show by induction on $n\in \NN$
that $G_n=0$. 
For $n=0$ we have $\sha_C^0(A)=A$.
So $G_0$ is in $A$. 
By Lemma~\ref{lem:fil} we have
\[ \Phi(\sum_{k\geq 1} G_k) \in F^1 \frakA(A).\]
Thus the first component of $\Phi(F)$ in $\frakA(A)$
is from
\[ \Phi(G_0)= G_0 \gamma_1 + {\rm\ a\ term\ in\ } F^1\frakA(A).\]
Thus $\Phi(G)=0$ implies $G_0\gamma_1=0$. 
Therefore, $G_0=0$. 

Now assume that $G_k=0$ for $k\leq n$ and consider
$G_{n+1}\in \sha_C^{n+1}(A)=A^{\otimes (n+2)}$.
Thus $G_{n+1}$ can be expressed as
\[ G_{n+1} = \sum_{i=1}^k a^{(i)}_1\otimes \ldots \otimes
    a^{(i)}_{n+2}, k\in \NN_+,\ a^{(i)}_j\in A. \]
Since $G_k=0$ for $k\leq n$, and by Lemma~\ref{lem:fil},
for $k\geq n+2$,

\begin{eqnarray*}
 \Phi(G_k)
&\in & \Phi(\sha^{k}_C(A))\\
&\subseteq & \Phi(\Fil^k \sha_C(A)) \\
&\subseteq & \Phi(\Fil^{n+2} \sha_C(A)) \\
&\subseteq & F^{n+2} \frakS(A), 
\end{eqnarray*}
the only contribution of $\Phi(G)$ to the coefficient of
$\gamma_{n+2}$ is from $\Phi (G_{n+1})$.
By Lemma~\ref{lem:lead}, this coefficient is
\[ \lambda^{n+1} \sum_{i=1}^k a^{(i)}_{n+2}\otimes \ldots \otimes
    a^{(i)}_1. \]
Thus $\Phi(F)=0$ implies that
\[ \lambda^{n+1}\sum_{i=1}^k a^{(i)}_{n+2}\otimes \ldots \otimes
    a^{(i)}_1 =0. \]
Since $\lambda$ is not a zero divisor in
$\bar{A}$, we further have 
\[ \sum_{i=1}^k a^{(i)}_{n+2}\otimes \ldots \otimes
    a^{(i)}_1 =0 \]
as an element in the tensor power algebra $A^{\otimes (n+2)}$.
But this element being zero or not depends only on the $C$-module
structure of $A^{\otimes (n+2)}$ and
the $C$-module map
\[ A^{\otimes (n+2)} \to A^{\otimes (n+2)},
a_{n+2}\otimes\ldots \otimes a_1 \mapsto
a_1 \otimes \ldots \otimes a_{n+2} \]
is an isomorphism.
Thus we also have
\[ G_{n+1} = \sum_{i=1}^k a^{(i)}_1\otimes \ldots \otimes
    a^{(i)}_{n+2} = 0. \]
This completes the induction.
Thus $\Phi$ is injective and hence an isomorphism. 
\proofend

\subsection{Completions}
We now give an internal construction of free complete Baxter
algebras by showing that the complete free Baxter algebra
$\widehat{\sha}_C(A)$ can be embedded into $\frakA(A)$.

Let $\frakA(A)\tp$ be the subgroup of $\frakA(A)$
consisting of sequences with finitely many non-zero entries.
Then we have
\[ \frakA(A)\tp = \oplus_{n=1}^\infty \overline{A} \gamma_n.\]
Define a filtration on $\frakA(A)\tp$ by taking 
\[F^k \frakA(A)\tp=\frakA(X) \cap F^k\frakA(A)\tp
 = \oplus_{n=k+1}^\infty \overline{A}\  \gamma_n. \]
By Proposition~\ref{prop:series}, we have
\[ \invlim (\frakA(A)\tp/\Fil^k \frakA(A)\tp)
    \cong \prod_{k\in\NN_+} \overline{A} \]
with the addition, multiplication and scalar multiplication defined
componentwise.
Therefore,
\[ \invlim (\frakA(A)\tp/\Fil^k \frakA(A)\tp) \cong \frakA(A).\]
From the definition of $F^k\frakA(A)$ and $F^k\frakA(A)\tp$, 
\[ \frakA(A)\tp/\Fil^n \frakA(A)\tp \cong
\frakA(A)/\Fil^n \frakA(A)\]
So $\frakA(A)$ is the completion of
itself with respect to the filtration $F^k \frakA(X)$. 
Since $F^k\frakS(A)=\frakS(A)\cap F^k\frakA(A)$, we have the
injective map of inverse systems
\[ \frakS(A)/F^k \frakS(A) \to \frakA(A)/F^k\frakA(A), k\in\NN_+.\]
So
\begin{equation}
\invlim (\frakS(A)/F^k\frakS(A)) \hookrightarrow 
\invlim (\frakA(A)/F^k\frakA(A)) \cong \frakA(A).
\label{eq:inv}
\end{equation}
We can easily describe the image of
$\invlim (\frakS(A)/F^k\frakS(A))$ in $\frakA(A)$.
It consists of sequences
$ (b^{(n)})_n,\ b^{(n)}\in \overline{A}$
that can be expressed as an infinite sum of the form
\[\sum_{k=1}^\infty (b^{(n)}_k)_n, \]
where $(b^{(n)}_k)_n\in \frakS(A)$ for each $k$. 
This means that, for any fixed $n\in\NN_+$,
all but finitely many $b^{(n)}_k,\ k\in \NN_+$,
are non-zero, and $\sum_{k=1}^\infty b^{(n)}_k=b^{(n)}$. 
We denote this image by $\tilde{\frakS}(A)$ with
the induced Baxter algebra structure. 

On the other hand, we also have the Baxter filtration
$\Fil^k$ on $\frakA(A)$ and $\frakS(A)$ (section~\ref{sec:comp}).

\begin{theorem}
\label{thm:rcomp}
\begin{enumerate}
\item
The Baxter algebra $\frakA(A)$ is complete.
\item
Assume that $\lambda\in C$ is not a zero divisor in
$\overline{A}$.
The isomorphism
$\Phi: \sha_C(A)\to \frakS(A)$
extends to an isomorphism of complete Baxter algebras
\[ \widehat{\Phi}: \widehat{\sha}_C(A) \to \tilde{\frakS}(A).\]
\end{enumerate}
\end{theorem}

We first prove a lemma. 
\begin{lemma}
\label{lem:fil2}
For any $k\in \NN_+$, we have
\begin{enumerate}
\item
$\Fil^k \frakA(A) = \lambda^ k F^k \frakA(A).$
\item
Assume that $\lambda\in C$ is not a zero divisor in
$\overline{A}$.
$\Fil^k \frakS(A) = F^k \frakS(A).$
\end{enumerate}
\end{lemma}

\proof
1.
We prove that, for any $n\in\NN_+$,
\[ P\tp_A(F^k\frakA(A))=\lambda F^{k+1}\frakA(A).\]
By the definition of $P\tp_A: \frakA(A)\to \frakA(A)$,
we have $P\tp_A(F^k\frakA(A))\subseteq \lambda F^{k+1} \frakA(A)$.
On the other hand, any element in $\lambda F^{k+1}\frakA(A)$
is of the form
\[ \lambda \sum_{i=k+2}^\infty a_i \gamma_i,\ a_i \in \overline{A}.\]
Then we have
\[ P\tp_A(\sum_{i=k+1}^\infty (a_{i+1}-a_i)\gamma_i)=
\lambda \sum_{i=k+2}^\infty a_i \gamma_i.\]
Here we take $a_{k+1}=0$.
This proves the equation. 

When $k=1$, we have $P\tp_A(\frakA(A)) =\lambda F^1\frakA(A)$.
Since $F^1\frakA(A)$ is already a Baxter ideal,
we have $\Fil^1 \frakA(A)=\lambda F^1\frakA(A)$. 
Inductively, assuming that
$\Fil^k \frakA(A) = \lambda^k F^k \frakA(A)$,
then we have
\[P\tp_A(\Fil^k \frakA(A)) = P\tp_A(\lambda^k F^k\frakA(A))
    =\lambda^{k+1} F^{k+1} \frakA(A).\]
Since $F^{k+1}\frakA(A)$ is a Baxter ideal of $\frakA(A)$, 
it is $\Fil^{k+1}\frakA(A)$, the Baxter ideal generated by
$P\tp_A(\Fil^k\frakA(A))$. 

2.
By Lemma~\ref{lem:fil}, Lemma~\ref{lem:lead} and
Theorem~\ref{thm:s-r}, we have,
for any $a\in \sha_C(A)$ 
\begin{eqnarray*}
&& \Phi(a) \in \Fil^k \frakS(A) \\
&\Leftrightarrow & a\in \Fil^k \sha_C(A) \\
&\Leftrightarrow & \Phi(a)\in \frakS(A)\cap F^k \frakA(A)\\
&\Leftrightarrow & \Phi(a)\in F^k \frakS(A).
\end{eqnarray*}
This proves the second equation.
\proofend

\noindent
{\bf Proof of Theorem~\ref{thm:rcomp}: }
From the first equation of Lemma~\ref{lem:fil2},
we have the exact sequence of inverse systems
\[ 0\to F^k \frakA(A/\lambda k A)\to \frakA(A)/\Fil^k \frakA(A)
\to \frakA(A)/F^k \frakA(A) \to 0.\]
Clearly
$\invlim F^k \frakA(A/\lambda k A) =0$.
Thus we have
\[ \hat{\frakA}(A))\hookrightarrow \invlim (\frakA(A)/F^k \frakA(A))
    =\frakA(A).\]
This proves the first statement.

Next assume that $\lambda\in C$ is not a zero divisor in
$\overline{A}$.
Then by the second statement of Lemma~\ref{lem:fil2},
\[ \hat{\frakS}(A) \cong \invlim (\frakS(A)/F^k \frakS(A)). \]
Then by Theorem~\ref{thm:s-r} and equation~(\ref{eq:inv})
we obtain
\[ \widehat{\sha}_C(A) \cong \tilde{\frakS}(A) \hookrightarrow
    \frakA(A). \]
\proofend

\addcontentsline{toc}{section}{\numberline {}References}

\end{document}